\documentclass[12pt,a4paper]{amsart}
\usepackage{graphicx}
\usepackage{epstopdf}
\usepackage{amssymb}
\usepackage{amsthm}
\usepackage{amsmath}
\usepackage{pstcol, pst-node}
\usepackage{mathrsfs}
\usepackage{amscd}

\usepackage{comment}

\usepackage[all]{xy}
\newtheorem{thm}{Theorem}

\newtheorem{cor}[thm]{Corollary}
\newtheorem{prop}[thm]{Proposition}
\newtheorem{lem}[thm]{Lemma}
\newtheorem{rem}[thm]{Remark}
\newtheorem{claim}[thm]{Claim}

\theoremstyle{definition}
\newtheorem{defn}[thm]{Definition}

\newtheorem{prop-def}[thm]{Proposition-Definition}

\newcommand{\disp}{\displaystyle}

\newcommand{\Ext}{\mathop{\mathrm{Ext}}\nolimits}

\newcommand{\Spec}{\mathop{\mathrm{Spec}}\nolimits}

\theoremstyle{remark}

\begin{document}
\author{Takeo Nishinou}
\date{}
\thanks{email : nishinou@rikkyo.ac.jp}
\address{Department of Mathematics, Rikkyo University,
	3-34-1, Nishi-Ikebukuro, Toshima, Tokyo, Japan} 
\subjclass[2000]{}
\keywords{}
\title{Deformation of pairs and semiregularity}
\maketitle
\begin{abstract}
In this paper, we study relative deformations of maps into a family of K\"ahler manifolds
 whose images are divisors.
We show that if the map satisfies a condition called semiregularity, 
 then it allows relative deformations if and only if the cycle class of the image remains Hodge
 in the family.
This gives a refinement of the so-called variational Hodge conjecture.
We also show that the semiregularity of maps is related to classical notions such as 
 Cayley-Bacharach conditions and d-semistability.
\end{abstract}
\section{Introduction}

Let $\pi\colon \mathfrak X\to D$
 be a deformation of a compact K\"ahler manifold $X_0$ of dimension $n\geq 2$
 over a disk $D$ in the complex plane.
Let $C_0$ be a compact reduced curve (when $n = 2$) or a compact smooth complex manifold
 of dimension $n-1$ (when $n>2$).
Let $\varphi_0\colon C_0\to X_0$ be a map which is an immersion,
 that is, for any $p\in C_0$, there is an open neighborhood $p\in V_p\subset C_0$
 such that $\varphi_0|_{V_p}$ is an embedding. 
Then the image of $\varphi_0$ determines an integral cohomology class $[\varphi_0(C_0)]$ of type $(1, 1)$, 
 that is, a Hodge class which is the Poincar\`e dual of the cycle $\varphi_0(C_0)$.
Note that the class $[\varphi_0(C_0)]$ naturally determines an integral cohomology class of each fiber of
 $\pi$.
Therefore, it makes sense to ask whether this class remains Hodge in these fibers or not.
Clearly, the condition that the class $[\varphi_0(C_0)]$ remains Hodge is necessary 
 for the existence of deformations of the map $\varphi_0$ to other fibers.

If we assume that the image is a local complete intersection, we can talk about the 
 \emph{semiregularity} of the map $\varphi_0$, see Section \ref{sec:2}.
When $\varphi_0$ is the inclusion, Bloch \cite{B} proved that if $\varphi_0$ is semiregular
 and the class $[\varphi_0(C_0)]$ remains Hodge, then there is a deformation of $\varphi_0$
 (Bloch proved it for local complete intersections of any codimension).
In other words, a local complete intersection subvariety which is semiregular
 satisfies the \emph{variational Hodge conjecture}. 
More precisely, the variational Hodge conjecture asks the existence of 
 a family of cycles of the class $[\varphi_0(C_0)]$ which need not restrict to $\varphi_0(C_0)$
 on the central fiber.
Therefore, Bloch's theorem in fact shows that the semiregularity gives a result stronger than 
 the variational Hodge conjecture.
However,  although Bloch's theorem guarantees the existence of a relative deformation of a cycle on the 
 central fiber $X_0$, it gives little control of the geometry of the deformed cycle.

Our purpose is to show that the semiregularity in fact suffices to control the geometry 
 of the deformed cycles when the cycle is of codimension one.
\begin{thm}\label{thm:1}
Assume that the map $\varphi_0$ is semiregular.
If the class $[\varphi_0(C_0)]$ remains Hodge, then the map $\varphi_0$ deforms
 to other fibers.
\end{thm}

For example, if the image $\varphi_0(C_0)$ has normal crossing singularity, 
 then there is a natural map $\widetilde\varphi_0\colon \widetilde C_0\to X_0$, 
 where $\widetilde C_0$ is the normalization of $C_0$ (when $n>2$, $C_0 = \widetilde C_0$).
Then if $\widetilde\varphi_0$ is semiregular, Theorem \ref{thm:1} implies that
 it deforms to a general fiber and the singularity of 
 the image remains the same (e.g., it gives a relative equigeneric deformation when $n = 2$).

On the other hand, if the image $\varphi_0(C_0)$ has normal crossing singularity, 
 the semiregularity turns out to be related to some classical notions appeared in different contexts. 
Namely, we will prove the following (see Corollary \ref{cor:3}).
\begin{thm}\label{thm:2}
Assume that the subvariety $\varphi_0(C_0)$ is semiregular in the classical sense.
That is, the inclusion of $\varphi_0(C_0)$ into $X_0$ is semiregular in the sense of Definition \ref{def:semiregular}.
Then if the map 
$H^{0}(\varphi_0(C_0), \mathcal N_{\iota}) \to H^{0}(\varphi_0(C_0), \mathcal S)$
 is surjective, the map $\varphi_0$ is semiregular.
In particular, if the class $[\varphi_0(C_0)]$ remains Hodge on the fibers of $\mathfrak X$, 
 the map $\varphi_0$ can be deformed to general fibers of $\mathfrak X$.
\end{thm}
Here $\mathcal S$ is the \emph{infinitesimal normal sheaf} of the variety $\varphi_0(C_0)$,
 see Section \ref{sec:6} for the definition.
A variety with normal crossing singularity is called \emph{d-semistable} 
 if the infinitesimal normal sheaf is trivial, see \cite{F}.
The notion of d-semistablity is known to be related to the existence of log-smooth deformations (see \cite{KF, KK})
By the above theorem, it turns out that it is also related to deformations of pairs, see Corollary \ref{cor:4}.\\

In the case where $n = 2$, Theorem 31 in \cite{N6} combined with Theorem \ref{thm:1} above
 implies the following.
Let $\varphi_0\colon C_0\to X_0$ be an immersion where $\varphi_0(C_0)$ is a reduced nodal curve.
Let $p\colon C_0\to \varphi_0(C_0)$ be the natural map (which is a partial normalization of $\varphi_0(C_0)$)
 and $P = \{p_i\}$ be the set of nodes of $\varphi_0(C_0)$ whose inverse image by $p$ consists of two points.
\begin{thm}\label{thm:CB}
	Assume that $\varphi_0(C_0)$ is semiregular in the classical sense
	and the class $[\varphi_0(C_0)]$ remains Hodge on the fibers of $\mathfrak X$.
	Then the map $\varphi_0$ deforms to general fibers of $\mathfrak X$
    if for each $p_i\in P$, there is a first order deformation of 
	$\varphi_0(C_0)$ which smoothes $p_i$, but does not smooth the other nodes of $P$. \qed
\end{thm}
This is related (in a sense opposite) to the classical \emph{Cayley-Bacharach condition}, see \cite{BHPV},
 which requires that
 if a first order deformation does not smooth the nodes $P\setminus \{p_i\}$, 
 then it does not smooth $p_i$, either.
Using this, we can also deduce a geometric criterion for the existence of deformations of pairs, 
 see Corollary \ref{cor:geomCB}.

Finally, based on a similar idea, we will prove that any projective variety can be swept by 
 nodal curves with very large number of nodes.
Namely, we will prove the following (see Corollary \ref{cor:20}).
\begin{thm}
Let $Y$ be a projective variety of dimension $n\geq 2$.
Then for any positive number $\varepsilon$, 
 there is an $(n-1)$-dimensional family $\mathcal C\to B$ of irreducible nodal curves
 whose fibers satisfy $\delta > g^{2-\varepsilon}$, and a map 
 $p\colon \mathcal C\to Y$ which dominates $Y$.
Here $\delta$ is the number of nodes of a fiber of $\mathcal C$ and 
 $g$ is the geometric genus of it.
\end{thm}
In general, it would be difficult to improve the exponent $2-\varepsilon$ further.
For example, if we can prove the existence of a nodal curve of large degree which satisfies the estimate 
 $\delta>g^{2+\varepsilon}$ when $Y$ is a Fano manifold, such a curve has many deformations 
 enough to carry out Mori's famous bend-and-break procedure \cite{Mo}.

\section{Semiregularity for local embeddings}\label{sec:2}

Let $n$ and $p$ be positive integers with $p<n$.
Let $M$ be a complex variety (not necessarily smooth or reduced) of dimension $n-p$
 and $X$ a compact K\"ahler manifold of dimension $n$.
Let $\varphi\colon M\to X$ be a map which is an immersion, 
 that is, for any $p\in M$, there is an open neighborhood $p\in U_p\subset M$
 such that $\varphi|_{U_p}$ is an embedding. 
We assume that the image is a local complete intersection.
Then, the normal sheaf $\mathcal N_{\varphi}$ is locally free of rank $p$.
Define the locally free sheaves $\mathcal K_{\varphi}$ and $\omega_M$ on $M$ by
\[
\mathcal K_{\varphi} = \wedge^p \mathcal N_{\varphi}^{\vee}
\]
and
\[
\omega_M = \mathcal K_{\varphi}^{\vee}\otimes \varphi^*\mathcal K_X,
\]
 where $\mathcal K_X$ is the canonical sheaf of $X$.

When $\varphi$ is an inclusion, the natural inclusion 
\[
\varepsilon\colon \mathcal N^{\vee}_{\varphi}\to \varphi^*\Omega_X^1
\]
gives rise to an element
\[
\begin{array}{ll}
\wedge^{p-1}\varepsilon \in Hom_{\mathcal O_M}\left(
\wedge^{p-1}\mathcal N_{\varphi}^{\vee}, \varphi^*\Omega_X^{p-1} \right)
&= \Gamma\left(M, (\varphi^*\Omega_X^{n-p+1})^{\vee}\otimes \varphi^*\mathcal K_X
 \otimes \mathcal K_{\varphi}^{\vee}\otimes\mathcal N_{\varphi}^{\vee}  \right)\\
&= Hom_{\mathcal O_X}(\Omega_X^{n-p+1}, \omega_M\otimes\mathcal N_{\varphi}^{\vee}).
\end{array}
\]

This induces a map on cohomology 
\[
\wedge^{p-1}\varepsilon\colon H^{n-p-1}(X, \Omega_X^{n-p+1})\to H^{n-p-1}(M, \omega_M\otimes \mathcal N_{\varphi}^{\vee}).
\]
When $\varphi$ is not an inclusion, then 
 $\Gamma\left(M, (\varphi^*\Omega_X^{n-p+1})^{\vee}\otimes \varphi^*\mathcal K_X
 \otimes \mathcal K_{\varphi}^{\vee}\otimes\mathcal N_{\varphi}^{\vee}  \right)$
 is not necessarily isomorphic to 
 $Hom_{\mathcal O_X}(\Omega_X^{n-p+1}, \omega_M\otimes\mathcal N_{\varphi}^{\vee})$, 
 but the map 
 $\wedge^{p-1}\varepsilon\colon H^{n-p-1}(X, \Omega_X^{n-p+1})\to 
 H^{n-p-1}(M, \omega_M\otimes \mathcal N_{\varphi}^{\vee})$
 is still defined.
\begin{defn}\label{def:semiregular}
We call $\varphi$ \emph{semiregular} if the natural map $\wedge^{p-1}\varepsilon$ is surjective.
\end{defn}

In this paper, we are interested in the case where $p = 1$ and $M$ is reduced when $n=2$, 
 and $M$ is smooth when $n>2$.
In this case, we have $\omega_M\otimes \mathcal N_{\varphi}^{\vee}\cong \varphi^*\mathcal K_X$ and the map 
 $\wedge^{p-1}\varepsilon$ will be
\[
 H^{n-2}(X, \mathcal K_X)\to H^{n-2}(M, \varphi^*\mathcal K_X).
\]

\section{Local calculation}
Let $\pi\colon \mathfrak X\to D$
 be a deformation of a compact K\"ahler manifold $X_0$ of dimension $n\geq 2$.
Here $D$ is a disk on the complex plane centered at the origin.
Let 
\[
\{(U_i, (x_{i, 1}, \dots, x_{i, n})\}
\]
 be a coordinate system of $X_0$.
Taking $D$ small enough, 
 the sets 
\[
(\mathfrak U_{i} = U_i\times D, (x_{i, 1}, \dots, x_{i, n}, t))
\]
 gives
 a coordinate system of $\mathfrak X$.
Precisely, we fix an isomorphism between 
$\mathfrak U_i$ and a suitable open subset of $\mathfrak X$ which is compatible with
$\pi$ and the inclusion $U_i\to X_0$.
Here $t$ is a coordinate on $D$ pulled back to $\mathfrak U_i$.
The functions $x_{i, l}$ are also pulled back to $\mathfrak U_i$ from $U_i$
by the natural projection. 

Take coordinate neighborhoods $\mathfrak U_i, \mathfrak U_j$ and $\mathfrak U_k$.
On the intersections of these open subsets, the coordinate functions on one of them can be written 
in terms of another.
We write this as follows.
Namely, on $\mathfrak U_i\cap \mathfrak U_j$, $x_{i, l}$ can be written as $x_{i, l}({\bf x}_{j}, t)$, 
 here we write 
\[
{\bf x}_j = (x_{j, 1}, \dots, x_{j, n}).
\]
Similarly, on $\mathfrak U_j\cap \mathfrak U_k$, we have $x_{j, l} = x_{j, l}({\bf x}_{k}, t)$.
Then, on $\mathfrak U_i\cap \mathfrak U_k$, we have
\[
x_{i, l} = x_{i, l}({\bf x}_{k}, t) = x_{i, l}({\bf x}_{j}({\bf x}_k, t), t).
\]

Let $X_t = \pi^{-1}(t)$ be the fiber of the family $\pi$ over $t\in D$.
Assume that there is a map 
\[
\varphi_0\colon C_0\to X_0
\]
 from a variety $C_0$ of dimension $n-1$ to $X_0$, which is an immersion.

We can take an open covering $\{V_{i}\}$ of $C_0$ such that 
 the restriction of $\varphi_0$ to $V_i$ is an embedding,
 the image $\varphi_0(V_i)$ is contained in $U_i$ and  
 is defined by an equation $f_{i, 0} = 0$ for some holomorphic function $f_{i, 0}$.
Let $\Spec\Bbb C[t]/t^{m+1}$ be the $m$-th order infinitesimal neighborhood fo the origin of $D$.
Note that 
\[
\{U_{i, m} = \mathfrak U_i\times_D \Spec\Bbb C[t]/t^{m+1}\}
\]
 gives a covering by coordinate neighborhoods
 of $X_m = \mathfrak X\times_D{\Spec\Bbb C[t]/t^{m+1}}$.  
We write by $x_{i, l, m}$ the restriction of $x_{i, l}$ to $U_{i, m}$.
Let us write 
\[
{\bf x}_{i, m} = \{x_{i, 1, m}, \dots, x_{i, n, m}\}.
\]

Assume we have constructed an $m$-th order deformation $\varphi_m\colon C_m\to \mathfrak X_m$
 of $\varphi_0$.
Here $m$ is a positive integer and
 $C_m$ is an $m$-th order deformation of $C_0$.
Let $V_{i, m}$ be the ringed space obtained by restricting $C_m$ to $V_i$.

Let $\{f_{i, m}({\bf x}_{i, m}, t)\}$
 be the set of local defining functions of $\varphi_m(V_{i, m})$ in $U_{i, m}$.
We will often write $f_{i, m}({\bf x}_{i, m}, t)$ as $f_{i, m}({\bf x}_{i, m})$ for notational simplicity.
In particular, on the intersection $U_{i, m}\cap U_{j, m}$, there is an invertible function $g_{ij, m}$ 
which satisfies 
\[
f_{i, m}({\bf x}_{i, m}({\bf x}_{j, m}, t), t) = g_{ij, m}({\bf x}_{j, m}, t)f_{j, m}({\bf x}_{j, m}, t) \;\; \text{(mod $t^{m+1}$)}.
\]



Define a holomorphic function $\nu_{ij, m}$ on $V_i\cap V_j$ by 
\[
t^{m+1}\nu_{ij, m}({\bf x}_{j, m+1}) = t^{m+1}\nu_{ij, m}({\bf x}_{j, 0}) = f_{i, m}({\bf x}_{i, m+1}({\bf x}_{j, m+1})) 
    - g_{ij, m}({\bf x}_{j, m+1})f_{j, m}({\bf x}_{j, m+1}),
\] 
which is an equality over $\Bbb C[t]/t^{m+2}$.
\begin{prop}
	The set of local sections $\{\nu_{ij,m}\}$ gives
	a \v{C}ech 1-cocycle with values in $\mathcal N_{\varphi_0}$ on $C_0$.
Here $\mathcal N_{\varphi_0}$ is the normal sheaf of the map $\varphi_0$ on $C_0$.  
\end{prop}
\begin{rem}
The cohomology class of the 
 cocycle $\{\nu_{ij, m}\}$ represents the obstruction to deforming the map $\varphi_m$ one step further.
\end{rem}
\proof
Note that $\mathcal N_{\varphi_0}$ is an invertible sheaf and the functions $g_{ij, 0}$ gives the
 transition functions of it.
Therefore, we need to check the identities
\[
\nu_{ik, m}({\bf x}_{k, m+1}) = \nu_{ij, m}({\bf x}_{j, m+1}({\bf x}_{k, m+1})) 
       + g_{ij, 0}({\bf x}_{j, 0}({\bf x}_{k, 0}))\nu_{jk, m}({\bf x}_{k, m+1})
\]
 and 
\[ 
\nu_{ij,m} = -g_{ij. 0}\nu_{ji, m}
\]
on $C_0$.

Note that 
\[
{\bf x}_{i, m+1}({\bf x}_{k, m+1}) \equiv {\bf x}_{i, m+1}({\bf x}_{j, m+1}({\bf x}_{k, m+1})) \;\; \text{mod} \;t^{m+2}.
\]
Then,
\[
\hspace{-.5in}\begin{array}{ll}
\disp t^{m+1}\nu_{ik, m}({\bf x}_{k, m+1})  & 
   \disp =  f_{i, m}({\bf x}_{i, m+1}({\bf x}_{k, m+1})) - g_{ik, m}({\bf x}_{k, m+1})f_{k, m}({\bf x}_{k, m+1})\\
& \disp = f_{i, m}({\bf x}_{i, m+1}({\bf x}_{j, m+1}({\bf x}_{k, m+1}))) 
   - g_{ij, m}({\bf x}_{j, m+1}({\bf x}_{k, m+1}))f_{j, m}({\bf x}_{j, m+1}({\bf x}_{k, m+1}))\\
& \hspace{.4in} + g_{ij, m}({\bf x}_{j, m+1}({\bf x}_{k, m+1}))f_{j, m}({\bf x}_{j, m+1}({\bf x}_{k, m+1}))
          - g_{ik, m}({\bf x}_{k, m+1})f_{k, m}({\bf x}_{k, m+1})\\
& \disp =    t^{m+1}\nu_{ij, m}({\bf x}_{j, m+1}({\bf x}_{k, m+1}))
                      + g_{ij, m}({\bf x}_{j, m+1}({\bf x}_{k, m+1}))(f_{j, m}({\bf x}_{j, m+1}({\bf x}_{k, m+1}))\\
                       & \hspace{.4in}   - g_{jk, m}({\bf x}_{k, m+1})f_{k, m}({\bf x}_{k, m+1}))
          \disp +g_{ij, m}({\bf x}_{j, m+1}({\bf x}_{k, m+1}))g_{jk, m}({\bf x}_{k, m+1})f_{k, m}({\bf x}_{k, m+1})\\
         &     \hspace{.4in}      - g_{ik, m}({\bf x}_{k, m+1})f_{k, m}({\bf x}_{k, m+1})\\
& \disp =  t^{m+1}\nu_{ij, m}({\bf x}_{j, m+1}({\bf x}_{k, m+1})) 
        + t^{m+1}g_{ij, m}({\bf x}_{j, m+1}({\bf x}_{k, m+1}))\nu_{jk, m}({\bf x}_{k, m+1})\\
       & \hspace{.4in} \disp + (g_{ij, m}({\bf x}_{j, m+1}({\bf x}_{k, m+1}))g_{jk, m}({\bf x}_{k, m+1})
             - g_{ik, m}({\bf x}_{k, m+1}))f_{k, m}({\bf x}_{k, m+1}).
\end{array}
\] 
Since
\[
(g_{ij, m}({\bf x}_{j, m+1}({\bf x}_{k, m+1}))g_{jk, m}({\bf x}_{k, m+1})
- g_{ik, m}({\bf x}_{k, m+1}))f_{k, m}({\bf x}_{k, m+1}) \equiv 0 \;\; \text{mod}\; t^{m+1}.
\]
 we have
\[
g_{ij, m}({\bf x}_{j, m+1}({\bf x}_{k, m+1}))g_{jk, m}({\bf x}_{k, m+1}) \equiv g_{ik, m}({\bf x}_{k, m+1})\;\; \text{mod} \;t^{m+1}
\]
Therefore, we have
\[
\begin{array}{l}
(g_{ij, m}({\bf x}_{j, m+1}({\bf x}_{k, m+1}))g_{jk, m}({\bf x}_{k, m+1})
- g_{ik, m}({\bf x}_{k, m+1}))f_{k, m}({\bf x}_{k, m+1}) \\
\hspace{.5in} \equiv (g_{ij, m}({\bf x}_{j, m+1}({\bf x}_{k, m+1}))g_{jk, m}({\bf x}_{k, m+1})
 - g_{ik, m}({\bf x}_{k, m+1}))f_{k, 0}({\bf x}_{k, m+1})  \;\; \text{mod}\; t^{m+2}.
\end{array}
\]
Since $f_{k, 0}(x_k) = 0$ on $C_0$, we have the first identity.
The second identity follows from this by taking $k = i$.\qed

\section{Explicit description of the Kodaira-Spencer class}
Let $\pi\colon \mathfrak X\to D$
 be a deformation of a smooth manifold $X_0$ as before.
We have the exact sequence
\[
0 \to \pi^*\Omega^1_{D}\to \Omega^1_{\mathfrak X}\to \Omega^1_{\mathfrak X/D}\to 0
\]
The Kodaira-Spencer class is, by definition, the corresponding class in 
 $\mu\in\Ext^1(\Omega^1_{\mathfrak X/D}, \pi^*\Omega^1_{D})$.


\begin{lem}
The class $\mu$ is represented by the \v{C}ech 1-cocycle
 $\mu_{ij} = \sum_{l=1}^n \frac{\partial x_{i, l}({\bf x}_j, t)}{\partial t}\partial_{x_{i, l}}dt$.
\end{lem}
\proof 
See \cite{Griffiths2}, Section II.1. \qed\\

From now on, we drop $dt$ from these expressions since it plays no role below. 
Restricting this to a presentation over $\Bbb C[t]/t^{m+1}$, we obtain the Kodaira-Spencer class for the 
 deformation $X_{m+1}:= \mathfrak X\times_D\Spec\Bbb C[t]/t^{m+2}$.
We denote this class by $\mu_m$.

Assume that we have constructed an $m$-th order deformation $\varphi_m\colon C_m\to X_m$.
Let $\mathcal N_{m/D}$ be the relative normal sheaf of $\varphi_m$
 and 
\[
p_m: \varphi_m^*\mathcal T_{X_{m}/D}\to \mathcal N_{m/D}
\]
 be the natural projection, 
 where $\mathcal T_{X_{m}/D}$ is the relative tangent sheaf of $X_m$.
Pulling $\mu_m$ back to $C_m$ and taking the image by $p_m$, we obtain a class
 $\overline{\mu}_m\in H^1(C_m, \mathcal N_{m/D})$.

As before, let $\{f_{i, m}({\bf x}_{i, m}, t)\}$
 be the set of local defining functions of $\varphi_m(V_{i, m})$ on $U_{i, m}$.
\begin{lem}
The class $\overline\mu_m$ is represented by the pull back of 
\[
\eta_{ij, m} =  \sum_{l=1}^n \frac{\partial x_{i, l}({\bf x}_j, t)}{\partial t}\partial_{x_{i, l}}f_{i, m}({\bf x}_i, t).
\]
 to $C_m$
\end{lem}
\proof
We check the cocycle condition.
Namely, we have
\[
\begin{array}{l}
\eta_{ik, m}-\eta_{ij, m}-g_{ij, m}\eta_{jk, m}  \\
 = \sum_{l=1}^n \frac{\partial x_{i, l}({\bf x}_k, t)}{\partial t}\partial_{x_{i, l}}f_{i, m}({\bf x}_i, t)
     - \sum_{l=1}^n \frac{\partial x_{i, l}({\bf x}_j, t)}{\partial t}\partial_{x_{i, l}}f_{i, m}({\bf x}_i, t)
   -g_{ij, m}\sum_{l=1}^n \frac{\partial x_{j, l}({\bf x}_k, t)}{\partial t}\partial_{x_{j, l}}f_{j, m}({\bf x}_i, t)\\
   = \sum_{l=1}^n \frac{\partial x_{i, l}({\bf x}_k, t)}{\partial t}\partial_{x_{i, l}}f_{i, m}({\bf x}_i, t)\\
    \hspace{.5in}- \sum_{l=1}^n \frac{\partial x_{i, l}({\bf x}_j, t)}{\partial t}\partial_{x_{i, l}}f_{i, m}({\bf x}_j, t)
    -g_{ij, m}\sum_{l=1}^n \frac{\partial x_{j, l}({\bf x}_k, t)}{\partial t}\partial_{x_{j, l}}(g_{ij, m}^{-1}f_{i, m}({\bf x}_i({\bf x}_j, t), t))\\
  = (\mu_{ik}-\mu_{ij}-\mu_{jk})f_{i, m} 
         -g_{ij, m}f_{i, m}({\bf x}_i({\bf x}_j, t), t)\sum_{l=1}^n \frac{\partial x_{j, l}({\bf x}_k, t)}{\partial t}\partial_{x_{j, l}}(g_{ij, m}^{-1}).
\end{array}
\]
Since $f_{i, m}({\bf x}_i({\bf x}_j, t), t)$ is zero on the image of $\varphi_m$, 
 we see that $\eta_{ik, m} = \eta_{ij, m} + g_{ij, m}\eta_{jk, m}$ on $C_m$.
Also, note that $g_{ij, m}$ is the transition function of the normal sheaf $\mathcal N_{m/D}$.
Then it is clear that $\eta_{ij, m}$ represents the class $\overline{\mu}_m$.\qed\\

Recall that an analytic cycle of codimension $r$ in a K\"ahler manifold
 determines a cohomology class of type $(r, r)$, which is the Poincar\'e dual of the 
 homology class of the cycle.
Let $\zeta_{C_0}\in H^{1}(X_0, \Omega_{X_0/D}^{1})$ be the class 
 corresponding to the image of $\varphi_0$.
Note that since the family $\mathfrak X$ is differential geometrically trivial, the class $\zeta_{C_0}$ 
 determines a cohomology class in $H^2(\mathfrak X, \Bbb C)$.
We denote it by $\widetilde\zeta_{C_0}$.
Then we have the following.
\begin{lem}
When $\varphi_0$ is semiregular, 
 the class $\widetilde\zeta_{C_0}$ remains Hodge in $X_{m+1}$ if and only if the class $\overline\mu_m$ is zero.
\end{lem}
\proof 
Since we are assuming we have constructed $\varphi_m\colon C_m\to X_m$, 
 the class $\widetilde\zeta_{C_0}$ is Hodge on $X_m$.
That is, $\widetilde\zeta_{C_0}|_{X_m}\in H^1(X_m, \Omega^1_{X_m/D})$.
In \cite[Proposition 4.2]{B}, Bloch shows that 
 $\widetilde\zeta_{C_0}$ remains Hodge on $X_{m+1}$ if and only if 
 the cup product $\widetilde\zeta_{C_0}|_{X_m}\cup \mu_m\in H^2(X_m, \mathcal O_{X_m})$ is zero.
This is the same as the claim that the cup product
 $\widetilde\zeta_{C_0}|_{X_m}\cup \mu_m\cup \alpha$ is zero for any $\alpha\in H^{2n-2}(X_m, \Bbb C)$.
On the other hand, we have the following.

\begin{claim}
The cup product 
 $\widetilde\zeta_{C_0}|_{X_m}\cup \mu_m\cup \alpha$ is zero for any $\alpha\in H^{2n-2}(X_m, \Bbb C)$
 if and only if
 the cup product $\overline\mu_m\cup \varphi_m^*\alpha$ is zero on $C_m$. 
\end{claim}
\noindent
{\it Proof of the claim.}
By the definition of $\widetilde\zeta_{C_0}|_{X_m}$, 
 the class $\widetilde\zeta_{C_0}|_{X_m}\cup \mu_m\cup \alpha$ is zero if and only if
 the class $\varphi_m^*\mu_m\cup \varphi_m^*\alpha$ is zero.
Note that the cohomology group $H^{2n-2}(X_m, \Bbb C)$ decomposes as
\[
H^{2n-2}(X_m, \Bbb C) \cong H^n(X_m, \Omega_{X_m/D}^{n-2})
    \oplus H^{n-1}(X_m, \Omega_{X_m/D}^{n-1})\oplus H^{n-2}(X_m, \mathcal K_{X_m/D}).
\]
By dimensional reason, the cup product between $\varphi_m^*\mu_m$
 and the pull back of the classes in $H^n(X_m, \Omega_{X_m/D}^{n-2})\oplus H^{n-1}(X_m, \Omega_{X_m/D}^{n-1})$
 is zero.
Therefore, we can assume that $\alpha$ belongs to $H^{n-2}(X_m, \mathcal K_{X_m/D})$, and so
 the class $\varphi_m^*\alpha$ belongs to $H^{n-2}(C_m, \varphi^*_m\mathcal K_{X_m/D})$.
On the other hand, $\varphi_m^*\mu_m$ belongs to 
 $H^1(C_m, \varphi_m^*\mathcal T_{X_m/D})$ and we have the natural map
\[
H^1(C_m, \varphi^*\mathcal T_{X_m/D})\to H^1(C_m, \mathcal N_{m/D}).
\]
Here $\overline{\mu}_m$ is the image of $\varphi_m^*\mu_m$ by this map.
Recall that the dual of 
 $H^1(C_m, \mathcal N_{m/D})$ is given by $H^{n-2}(C_m, \varphi_m^*\mathcal K_{X_m/D})$.
Therefore, it follows that the cup product $\varphi^*\mu_m\cup \varphi^*\alpha$
 reduces to $\overline\mu_m\cup \varphi_m^*\alpha$.
This proves the claim.\qed\\

It immediately follows that if $\overline{\mu}_m$ is zero, then 
 $\widetilde\zeta_{C_0}$ remains Hodge in $X_{m+1}$.
For the converse, assume that $\widetilde\zeta_{C_0}$ remains Hodge in $X_{m+1}$.
%
There is a natural map 
\[
\iota: H^{2n-2}(X_m, \Bbb C) \to H^1(C_m, \mathcal N_{m/D})^{\vee}
\]
 as in the proof of the claim.
Namely, 
 for a class $\alpha$ of 
 $H^{2n-2}(X_m, \Bbb C)=
 H^{n}(X_m, \Omega_{X_m/D}^{n-2})\oplus H^{n-1}(\Omega_{X_m/D}^{n-1})
   \oplus H^{n-2}(\Omega_{X_m/D}^{n})$
   and $\beta\in H^1(C_m, \mathcal N_{m/D})$, 
 let
 $\iota(\alpha)(\beta)$ be the cup product $\beta\cup \varphi_m^*\alpha$
 composed with the trace map $H^{n-1}(C_m, \omega_{C_m})\to \Bbb C$.

The restriction of this map to $X_0$ is a surjection by the semiregularity of $\varphi_0$.
Since the surjectivity is an open condition, $\iota$ is also a surjection.
This shows that $\overline\mu_m\cup \varphi_m^*\alpha$ is zero
 for any $\alpha\in H^{2n-2}(X_m, \Bbb C)$ is equivalent to the claim that $\overline\mu_m$ is zero.\qed\\

Thus, in this case we can write $\overline{\mu}_m$ as the coboundary of a
 \v{C}ech 0-cochain with values in $\mathcal N_{m/D}$ on $C_m$.
We choose one such representative $\{\delta_i\}$ where 
 $\delta_i\in\Gamma(V_{i, m}, \mathcal N_{m/D})$. 
Also note that by the exact sequence 
\[
0\to \mathcal O_{U_{i, m}}\to \mathcal O_{U_{i, m}}(\varphi_m(V_{i, m}))\to \mathcal N_{m/D}|_{V_{i, m}}\to 0,  
\]
 there is a section $\widetilde\delta_i$ of $\mathcal O_{U_{i, m}}(\varphi_m(V_{i, m}))$
 which maps to $\delta_i$.
Explicitly, putting 
 $\widetilde\eta_{ij, m} =
    \widetilde\delta_i({\bf x}_i({\bf x}_j, t), t) - g_{ij, m}({\bf x}_j, t)\widetilde\delta_j({\bf x}_j, t)$, 
 it coincides with $\eta_{ij, m}$ when restricted to $V_{i, m}$.

\section{Proof of the theorem}
Recall that the obstruction to deforming $\varphi_m$ is given by a cocycle $\nu_{ij, m}$ on $C_0$
 defined by 
 $t^{m+1}\nu_{ij, m}({\bf x}_{j}) = f_{i, m}({\bf x}_{i}({\bf x}_{j}, t), t) 
 - g_{ij, m}({\bf x}_{j}, t)f_{j, m}({\bf x}_{j}, t)$.
Differentiating this with respect to $t$, we have
\[
\begin{array}{ll}
(m+1)t^m\nu_{ij, m}({\bf x}_{j}) & = \frac{\partial f_{i, m}({\bf x}_i, t)}{\partial t} 
 + \sum_{l=1}^n\frac{\partial x_{i, l}({\bf x}_j, t)}{\partial t}\frac{\partial f_{i, m}({\bf x}_i, t)}{\partial {x}_{i, l}}
  - g_{ij, m}({\bf x}_{j}, t)\frac{\partial f_{j, m}({\bf x}_j, t)}{\partial t}
  -\frac{\partial g_{ij, m}({\bf x}_j, t)}{\partial t}f_{j, m}({\bf x}_j, t)\\
 &= \frac{\partial f_{i, m}({\bf x}_i, t)}{\partial t} 
 - g_{ij, m}({\bf x}_{j}, t)\frac{\partial f_{j, m}({\bf x}_j, t)}{\partial t}
  + \eta_{ij, m}
 -\frac{\partial g_{ij, m}({\bf x}_j, t)}{\partial t}f_{j, m}({\bf x}_j, t).
\end{array}
\]
Since $f_{j, m}$ is zero on $C_m$, we can ignore the last term.
By the same reason, we can replace $\eta_{ij, m}$ by $\widetilde\eta_{ij, m}$.

Dividing this by $f_{i, m}({\bf x}_i, t)$, we have
\[
\begin{array}{ll}
(m+1)t^m\frac{\nu_{ij, m}({\bf x}_j)}{f_{i, m}({\bf x}_i, t)}
 & = \frac{1}{f_{i, m}({\bf x}_i, t)}\frac{\partial f_{i, m}({\bf x}_i, t)}{\partial t} 
         -  \frac{g_{ij, m}({\bf x}_{j}, t)}{f_{i, m}({\bf x}_i, t)}\frac{\partial f_{j, m}({\bf x}_j, t)}{\partial t}
         + \frac{\eta_{ij, m}}{f_{i, m}({\bf x}_i, t)}\\
   & = \frac{1}{f_{i, m}({\bf x}_i, t)}\frac{\partial f_{i, m}({\bf x}_i, t)}{\partial t} 
   -  \frac{g_{ij, m}({\bf x}_{j}, t)}{f_{i, m}({\bf x}_i, t)}\frac{\partial f_{j, m}({\bf x}_j, t)}{\partial t}
              +\frac{\widetilde\delta_{i}}{f_{i, m}({\bf x}_i, t)}-\frac{g_{ij, m}\widetilde\delta_{j}}{f_{i, m}({\bf x}_i, t)}\\
   & = \frac{1}{f_{i, m}({\bf x}_i, t)}(\frac{\partial f_{i, m}({\bf x}_i, t)}{\partial t} +\widetilde\delta_i)
          - \frac{1}{f_{j, m}({\bf x}_i, t)}(\frac{\partial f_{j, m}({\bf x}_j, t)}{\partial t} +\widetilde\delta_j)
\end{array}
\]
 modulo functions holomorphic on $C_m$.
Note that this is an equation over $\Bbb C[t]/t^{m+1}$, and so we have
 $\frac{g_{ij, m}({\bf x}_{j}, t)f_{j, m}({\bf x}_j, t)}{f_{i, m}({\bf x}_i, t)} = 1$.
Let $[\frac{1}{f_{i, m}({\bf x}_i, t)}(\frac{\partial f_{i, m}({\bf x}_i, t)}{\partial t} +\widetilde\delta_i)]_m$
 be the coefficient of $t^m$ in 
 $\frac{1}{f_{i, m}({\bf x}_i, t)}(\frac{\partial f_{i, m}({\bf x}_i, t)}{\partial t} +\widetilde\delta_i)$.
Note that the above equation still holds when we replace 
 $\frac{1}{f_{i, m}({\bf x}_i, t)}(\frac{\partial f_{i, m}({\bf x}_i, t)}{\partial t} +\widetilde\delta_i)$
 and $\frac{1}{f_{j, m}({\bf x}_i, t)}(\frac{\partial f_{j, m}({\bf x}_j, t)}{\partial t} +\widetilde\delta_j)$
 by $[\frac{1}{f_{i, m}({\bf x}_i, t)}(\frac{\partial f_{i, m}({\bf x}_i, t)}{\partial t} +\widetilde\delta_i)]_m$
 and $[\frac{1}{f_{j, m}({\bf x}_i, t)}(\frac{\partial f_{j, m}({\bf x}_j, t)}{\partial t} +\widetilde\delta_j)]_m$, 
 respectively.
Also, we can think $[\frac{1}{f_{i, m}({\bf x}_i, t)}(\frac{\partial f_{i, m}({\bf x}_i, t)}{\partial t} +\widetilde\delta_i)]_m$
 as a function on $U_i$ by forgetting $t^m$.


Now we proceed as in \cite{N6} to produce a family of 
 local $C^{\infty}$ differential forms on $C_0$ which represent the obstruction class $\{\nu_{ij, m}\}$. 
Namely, introduce any Hermitian metric on $X_0$.
For each $V_{i}$, the unit normal bundle of the image $\varphi_0(V_{i})$
 of radius $r$ gives a circle bundle on it.
Here $r$ is a small positive real number and when $C_0$ is a singular curve, then we 
 construct the bundle only away from a small neighborhood of singular points.
These local circle bundles glue and give a global circle bundle on $C_0$ (away from 
 singular points).

Note that the class $[\nu_{ij, m}]$ is zero if and only if the pairing of it with 
 any class in $H^{n-2}(C_0, \varphi^*\mathcal K_{X_0})$ is zero.
By semiregularity, any class in $H^{n-2}(C_0, \varphi^*K_{X_0})$
 is a restriction of an element of $H^{n-2}(X_0, \mathcal K_{X_0})$.
Let $\Theta$ be any closed $C^{\infty}$ $(2n-2)$-form on $X_0$.
In particular, $\Theta$ represents a class in 
\[
H^{2n-2}(X_0) = H^{n-2}(X_0, \mathcal K_{X_0})\oplus
  H^{n-1}(X_0, \Omega^{n-1}_{X_0})\oplus H^n(X_0, \Omega^{n-2}_{X_0}).
\]
Here $\Omega^i_{X_0}$ is the sheaf of holomorphic $i$-forms on $X_0$.
Integrating the restriction of the singular $(2n-2)$-form
 $[\frac{1}{f_{i, m}({\bf x}_i, t)}(\frac{\partial f_{i, m}({\bf x}_i, t)}{\partial t} +\widetilde\delta_i)]_m\Theta$ to the circle bundle
 along the fibers, we obtain a set of local closed $(2n-3)$-forms $\gamma_i$ on $V_{i}$.
As the radius $r$ goes to zero, the \v{C}ech 1-cocycle obtained as the differences of
 $\{\gamma_i\}$ converges to the obstruction class $[\nu_{ij, m}]$ paired with the pull back of $\Theta$
 by $\varphi_0$.

The argument in \cite{N6} proves the following:

\begin{enumerate}
\item The class obtained from $\{\gamma_i\}$ does not depend on the radius $r$.
\item Thus, if $C_0$ is nonsingular, then since $\{\gamma_i\}$ is defined on an open covering 
 of $C_0$, the class determined by $\{\gamma_i\}$ is zero by definition for any $\Theta$.
This implies that $[\nu_{ij, m}]\in H^1(C_0, \mathcal N_{\varphi_0})$ is zero, too.
\item If $C_0$ has singular points, $\{\gamma_i\}$ is defined only on away from the singular points, 
 so one cannot immediately conclude that the class $[\nu_{ij, m}]$ is zero.
However, we can identify the class determined by $\{\gamma_i\}$ by local calculation at singular points, 
 and Stokes' theorem shows each of these local contributions is zero, which implies 
 the class $[\nu_{ij, m}]\in H^1(C_0, \mathcal N_{\varphi_0})$ is again zero.
\end{enumerate}
This finishes the proof of the theorem.\qed

\section{Criterion for semiregularity}\label{sec:6}
In this section, we give necessary conditions for a map $\varphi_0\colon C_0\to X_0$
 to be semiregular.
It turns out that some classical notions which appeared in different contexts such as Cayley-Bacharach condition
 and d-semistability are related to relative deformations of maps. 
\subsection{The case $n>2$}
First we consider the case $n>2$.
Let $\pi\colon \mathfrak X\to D$ be a family of $n$-dimensional K\"ahler manifolds.
Let $\varphi_0\colon C_0\to X_0$ be a map from a compact smooth complex manifold of dimension 
 $n-1$ which is an immersion.
We also assume that the image $\varphi_0(C_0)$ has normal crossing singularity.
 
Consider the exact sequence on $\varphi_0(C_0)$ given by
\[
0\to \iota^*\mathcal K_{X_0}\to p_*\varphi_0^*\mathcal K_{X_0}\to \mathcal Q\to 0,
\] 
 where $\iota\colon \varphi_0(C_0)\to X_0$ is the inclusion, and $p\colon C_0\to \varphi_0(C_0)$
 is the normalization.
The sheaf $\mathcal Q$ is defined by this sequence.
It is supported on the singular locus $sing(\varphi_0(C_0))$ of $\varphi_0(C_0)$.
We have an associated exact sequence of cohomology groups
\[
\begin{array}{ll}
\cdots\to H^{n-2}(\varphi_0(C_0), \iota^*\mathcal K_{X_0})& \to H^{n-2}(\varphi_0(C_0), p_*\varphi_0^*\mathcal K_{X_0})
 \to H^{n-2}(\varphi_0(C_0), \mathcal Q)\\
 &  \to
 H^{n-1}(\varphi_0(C_0), \iota^*\mathcal K_{X_0})\to H^{n-1}(\varphi_0(C_0), p_*\varphi_0^*\mathcal K_{X_0})
 \to H^{n-1}(\varphi_0(C_0), \mathcal Q).
\end{array}
\]
By dimensional reason, we have $H^{n-1}(\varphi_0(C_0), \mathcal Q) = 0$.
Also, note that 
\[
H^{i}(\varphi_0(C_0), p_*\varphi_0^*\mathcal K_{X_0})\cong H^{i}(C_0, \varphi_0^*\mathcal K_{X_0})
\]
 for $i = n-2, n-1$, by the Leray's spectral sequence.
Therefore, 
 if $\varphi_0(C_0)$ is semiregular in the classical sense, that is, 
 the natural map $H^{n-2}(X_0, \mathcal K_{X_0})\to H^{n-2}(\varphi_0(C_0), \iota^*\mathcal K_{X_0})$
 is surjective, then the map $\varphi_0$ is semiregular if and only if the map
 $H^{n-2}(\varphi_0(C_0), \iota^*\mathcal K_{X_0}) \to H^{n-2}(\varphi_0(C_0), p_*\varphi_0^*\mathcal K_{X_0})$
 is surjective.
\begin{cor}\label{cor:1}
Assume that $\varphi_0(C_0)$ is semiregular in the classical sense
 and the class $[\varphi_0(C_0)]$ remains Hodge on the fibers of $\mathfrak X$.
Then if the map 
 $H^{n-2}(\varphi_0(C_0), \iota^*\mathcal K_{X_0}) \to H^{n-2}(\varphi_0(C_0), p_*\varphi_0^*\mathcal K_{X_0})$
 is surjective, $\varphi_0$ can be deformed to general fibers of $\mathfrak X$.\qed
\end{cor}

On the other hand, consider the exact sequence
\[
0\to p_*\mathcal N_{\varphi_0}\to \mathcal N_{\iota}\to \mathcal S\to 0,
\]
 of sheaves on $\varphi_0(C_0)$, where $\mathcal S$ is defined by this sequence.
The associated exact sequence of cohomology groups is
\[
\begin{array}{l}
0\to H^0(\varphi_0(C_0), p_*\mathcal N_{\varphi_0})\to H^0(\varphi_0(C_0), \mathcal N_{\iota})
 \to H^0(\varphi_0(C_0), \mathcal S)\\
   \hspace{.4in}  \to H^1(\varphi_0(C_0), p_*\mathcal N_{\varphi_0})\to H^1(\varphi_0(C_0), \mathcal N_{\iota})\to \cdots
\end{array}
\]
We have 
 $H^i(\varphi_0(C_0), p_*\mathcal N_{\varphi_0})\cong H^i(C_0, \mathcal N_{\varphi_0})$
 again by the Leray's spectral sequence.
Note that the group $H^i(C_0, \mathcal N_{\varphi_0})$
 is isomorphic to the dual of $H^{n-1-i}(C_0, \varphi_0^*\mathcal K_{X_0})$, $i = 0, 1$.
Similarly, the group $H^i(\varphi_0(C_0), \mathcal N_{\iota})$ is isomorphic to the dual of 
 $H^{n-1-i}(\varphi_0(C_0), \iota^*\mathcal K_{X_0})$, $i = 0, 1$.

\begin{comment}
(Since $\varphi_0(C_0)$ is reduced).
\end{comment}

Comparing the dual of the previous cohomology exact sequence with the latter, 
 we obtain $H^{n-2}(\varphi_0(C_0), \mathcal Q)^{\vee}\cong H^0(\varphi_0(C_0), \mathcal S)$.
In particular, we can restate Corollary \ref{cor:1} as follows.
\begin{cor}\label{cor:3}
Assume that $\varphi_0(C_0)$ is semiregular in the classical sense
and the class $[\varphi_0(C_0)]$ remains Hodge on the fibers of $\mathfrak X$.
Then if the map 
$H^{0}(\varphi_0(C_0), \mathcal N_{\iota}) \to H^{0}(\varphi_0(C_0), \mathcal S)$
 is surjective, $\varphi_0$ can be deformed to general fibers of $\mathfrak X$.\qed
\end{cor}

The sheaf $\mathcal S$ is the \emph{infinitesimal normal sheaf} of the singular locus of $\varphi_0(C_0)$, 
 as we will see below.
Recall that we assume that the image $\varphi_0(C_0)$ has normal crossing singularity.
Then, for any point $p\in \varphi_0(C_0)$, we can take a coordinate system 
 $(x_1, \dots, x_n)$ on a neighborhood $U$ of 
 $p$ in $X_0$ 
 so that $U\cap \varphi_0(C_0)$ is given by $x_1\cdots x_k = 0$, $1\leq k\leq n$.
Let $\mathcal I_j$ the ideal of $\mathcal O_U$ generated by $x_j$
 and let $\mathcal I$ be the ideal defining $\varphi_0(C_0)\cap U$ in $U$. 
Then 
\[
\mathcal I_1/\mathcal I_1\mathcal I\otimes \cdots \otimes \mathcal I_k/\mathcal I_k\mathcal I
\]
 gives an invertible sheaf on the singular locus of $\varphi_0(C_0)\cap U$.
Globalizing this construction, we obtain an invertible sheaf on the singular locus of $\varphi_0(C_0)$.
Then the dual invertible sheaf of this is called the infinitesimal normal sheaf of
 the singular locus of $\varphi_0(C_0)$, see \cite{F}.
 
\begin{lem}
The sheaf $\mathcal S$ is isomorphic to the infinitesimal normal sheaf.
\end{lem}
\proof
Note that the sheaf
 $\mathcal I_1/\mathcal I_1\mathcal I\otimes \cdots \otimes \mathcal I_k/\mathcal I_k\mathcal I$
 is generated by the element $x_1\otimes \cdots \otimes x_k$.
The sheaf $p_*\mathcal N_{\varphi_0}$ is given by 
 $\oplus_{i=1}^kHom(\mathcal I_i/\mathcal I_i^2, \mathcal O_U)$ on $U$.
The sheaf $\mathcal N_{\iota}$ is given by $Hom(\mathcal I/\mathcal I^2, \mathcal O_U)$.
The sheaf $\mathcal N_{\iota}$ is an invertible sheaf and generated by the morphism which 
 maps $x_1\cdots x_k$ to $1\in\mathcal O_U$.
In particular, by multiplying any $x_1\cdots \check{x}_i\cdots x_k$, the generator is mapped into the image of 
 $p_*\mathcal N_{\varphi_0}\to \mathcal N_{\iota}$, 
 namely, the image of the generator of $Hom(\mathcal I_i/\mathcal I_i^2, \mathcal O_U)$.
Also note that the ideal of the singular locus of $\varphi_0(C_0)$ is generated by 
 $x_1\cdots \check{x}_i\cdots x_k$, $i = 1, \dots, k$.
From these, it is easy to see that the cokernel of the map 
 $p_*\mathcal N_{\varphi_0}\to \mathcal N_{\iota}$ is isomorphic to the dual of 
 $\mathcal I_1/\mathcal I_1\mathcal I\otimes \cdots \otimes \mathcal I_k/\mathcal I_k\mathcal I$.\qed\\

Recall that the infinitesimal normal sheaf is related to deformations of $\varphi_0(C_0)$
 which smooth the singular locus,  see \cite{F}.
In particular, $\varphi_0(C_0)$ is called \emph{d-semistable} if the infinitesimal normal sheaf is trivial, 
 and d-semistable variety carries a log structure log smooth over a standard log point, so that 
 one can study its deformations via log smooth deformation theory \cite{KF, KK, KN}. 

By Corollary \ref{cor:3}, the infinitesimal normal sheaf plays a crucial in the deformation theory
 even if it is not d-semistable.
On the other hand, the notion of d-semistability gives a sufficient condition for the existence of 
 deformations in this situation, too, as follows.
\begin{cor}\label{cor:4}
Assume that the image $\varphi_0(C_0)$ is very ample and $H^1(X_0, \mathcal O_{X_0}(\varphi_0(C_0))) = 0$.
Assume also that $\varphi_0(C_0)$ is d-semistable and the singular locus of $\varphi_0(C_0)$ is connected.
Then the map $\varphi_0$ is semiregular.
\end{cor}
\proof
First, the subvariety $\varphi_0(C_0)$ is semiregular in the classical sense. 
Namely, consider the cohomology exact sequence
\[
\cdots \to H^1(X_0, \mathcal O_{X_0}(\varphi_0(C_0)))\to H^1(\varphi_0(C_0), \mathcal N_{\iota})
 \to H^2(X_0, \mathcal O_{X_0})\to \cdots, 
\]
 here $\iota\colon \varphi_0(C_0)\to X_0$ is the inclusion.
When  $H^1(X_0, \mathcal O_{X_0}(\varphi_0(C_0))) = 0$,
 the map $H^1(\varphi_0(C_0), \mathcal N_{\iota})\to H^2(X_0, \mathcal O_{X_0})$ is injective.
Since this map is the dual of the semiregularity map $H^{n-2}(X_0, \mathcal K_{X_0})\to 
 H^{n-2}(\varphi_0(C_0), \mathcal \iota^*\mathcal K_{X_0})$, 
 it follows that $\varphi_0(C_0)$ is semiregular.

To prove that $\varphi_0$ is semiregular, it suffices to show the map 
　$H^{0}(\varphi_0(C_0), \mathcal N_{\iota}) \to H^{0}(\varphi_0(C_0), \mathcal S)$
　is surjective.
When $\varphi_0(C_0)$ is d-semistable, the sheaf $\mathcal S$ is the trivial line bundle on the 
 singular locus of $\varphi_0(C_0)$. 
Since we assume that the singular locus is connected, it suffices to show that the map 
 $H^{0}(\varphi_0(C_0), \mathcal N_{\iota}) \to H^{0}(\varphi_0(C_0), \mathcal S)$ is not the zero map.
This in turn is equivalent to the claim that the injection 
 $H^0(\varphi_0(C_0), p_*\mathcal N_{\varphi_0})\to H^0(\varphi_0(C_0), \mathcal N_{\iota})$
 is not an isomorphism.
Since $\varphi_0(C_0)$ is very ample, there is a section $s$ of $\mathcal O_X(\varphi_0(C_0))$
 which does not entirely vanish on the singular locus of $\varphi_0(C_0)$.
Then if $\sigma$ is a section of $\mathcal O_X(\varphi_0(C_0))$ defining $\varphi_0(C_0)$, 
 the sections $\sigma + \tau s$, where $\tau\in\Bbb C$ is a parameter, deforms
 $\varphi_0(C_0)$, and the non-vanishing of $s$ on the singular locus of $\varphi_0(C_0)$ implies that
 this smoothes a part of the singular locus of $\varphi_0(C_0)$. 
Since the sections of $H^0(\varphi_0(C_0), p_*\mathcal N_{\varphi_0})$ give first order deformations
 which does not smooth the singular locus, it follows that the map 
 $H^0(\varphi_0(C_0), p_*\mathcal N_{\varphi_0})\to H^0(\varphi_0(C_0), \mathcal N_{\iota})$
 is not an isomorphism.
This proves the claim.\qed


\subsection{The case $n  = 2$}
Now let us consider the case $n = 2$.
Although we can work in a more general situation,
 we assume $\varphi_0(C_0)$ is a reduced nodal curve for simplicity.
However $C_0$ need not be smooth.
Let $p\colon C_0\to \varphi_0(C_0)$ be the natural map, which is a partial normalization.
In this case, we can deduce very explicit criterion for the semiregularity.
Again, we have the exact sequence
\[
\begin{array}{l}
0\to H^0(\varphi_0(C_0), p_*\mathcal N_{\varphi_0})\to H^0(\varphi_0(C_0), \mathcal N_{\iota})
\to H^0(\varphi_0(C_0), \mathcal S)\\
\hspace{.4in}  \to H^1(\varphi_0(C_0), p_*\mathcal N_{\varphi_0})\to H^1(\varphi_0(C_0), \mathcal N_{\iota})\to \cdots,
\end{array}
\]
 and if $\varphi_0(C_0)$ is semiregular in the classical sense,
 then $\varphi_0$ is semiregular if and only if the map 
 $H^0(\varphi_0(C_0), \mathcal N_{\iota})
 \to H^0(\varphi_0(C_0), \mathcal S)$ is surjective.
Let $P = \{p_i\}$ be the set of nodes of $\varphi_0(C_0)$ whose inverse image by $p$ consists of two points.
Then the sheaf $\mathcal S$ is isomorphic to $\oplus_{i}\Bbb C_{p_i}$, where
 $\Bbb C_{p_i}$ is the 
 skyscraper sheaf at $p_i$.
By an argument similar to the one in the previous subsection, we proved the following
  in \cite{N6}.
\begin{thm}\label{thm:CB}
	Assume that $\varphi_0(C_0)$ is semiregular in the classical sense.
	Then the map $\varphi_0$ is semiregular if and only if for each $p_i\in P$, there is a first order deformation of 
	$\varphi_0(C_0)$ which smoothes $p_i$, but does not smooth the other nodes of $P$. \qed
\end{thm}

For applications, it will be convenient to write this in a geometric form.
Consider the exact sequence
\[
0\to \mathcal O_{X_0}\to \mathcal O_{X_0}(\varphi_0(C_0))\to \mathcal N_{\iota}\to 0
\]
 of sheaves on $X_0$ and the associated cohomology sequence
\[
0 \to H^0(X_0, \mathcal O_{X_0})\to H^0(X_0, \mathcal O_{X_0}(\varphi_0(C_0)))\to 
 H^0(\varphi_0(C_0), \mathcal N_{\iota})
 \to H^1(X_0, \mathcal O_{X_0})\to \cdots.
\]
Let $V$ be the image of the map $H^0(\varphi(C_0), \mathcal N_{\iota})
 \to H^1(X_0, \mathcal O_{X_0})$.
Since we are working in the analytic category, we have the exact sequence
\[
0\to \Bbb Z\to \mathcal O_{X_0}\to \mathcal O_{X_0}^{\ast}\to 0
\]
 of sheaves on $X$.
Let $\bar V$ be the image of $V$ in $Pic^0(X_0) = H^1(X_0, \mathcal O_{X_0}^{\ast})$.
In \cite{N6}, we proved the following.

\begin{cor}\label{cor:geomCB}
	In the situation of Theorem \ref{thm:CB}, the map $\varphi_0$ is unobstructed
	if for each $p_i\in P$, there is an effective divisor $D$ such that $\mathcal O_X(\varphi_0(C_0)-D)\in \bar V$
	which avoids $p_i$ but passes through all points in $P\setminus\{p_i\}$. 
\end{cor}

A particularly nice case is
 when the map $H^0(\varphi_0(C_0), \mathcal N_{\iota})\to H^1(X_0, \mathcal O_{X_0})$ is surjective.
This is the case when $\varphi_0(C_0))$ is sufficiently ample.
Then, if for each $p_i\in P$ there is an effective divisor $D$ which is algebraically equivalent to $\varphi_0(C_0)$
 which avoids $p_i$ but passes through all points in $P\setminus\{p_i\}$, the map $\varphi_0$ is semiregular.
This is, in a sense, the opposite to the classical \emph{Cayley-Bacharach property},
 see for example \cite{BHPV}.

Combined with Theorem \ref{thm:1}, we have the following.
\begin{cor}\label{cor:5}
	Assume that $\varphi_0(C_0)$ is reduced, nodal and semiregular in the classical sense
	and the class $[\varphi_0(C_0)]$ remains Hodge on the fibers of $\mathfrak X$.
	Then the map $\varphi_0$ deforms to general fibers of $\mathfrak X$
	if the condition in Theorem \ref{thm:CB} or Corollary \ref{cor:geomCB} is satisfied. \qed
\end{cor}

In the case of $n = 2$, the original exact sequence 
\[
\begin{array}{ll}
\cdots\to H^{0}(\varphi_0(C_0), \iota^*\mathcal K_{X_0})& \to H^{0}(\varphi_0(C_0), p_*\varphi_0^*\mathcal K_{X_0})
\to H^{0}(\varphi_0(C_0), \mathcal Q)\\
&  \to
H^{1}(\varphi_0(C_0), \iota^*\mathcal K_{X_0})\to H^{1}(\varphi_0(C_0), p_*\varphi_0^*\mathcal K_{X_0})
\to H^{1}(\varphi_0(C_0), \mathcal Q)
\end{array}
\]
 before taking the dual is sometimes also useful.
In this case, if $\varphi_0(C_0)$ is semiregular in the classical sense, then $\varphi_0$ is semiregular
 if and only if the map 
\[
H^{0}(\varphi_0(C_0), \iota^*\mathcal K_{X_0}) \to H^{0}(\varphi_0(C_0), p_*\varphi_0^*\mathcal K_{X_0})
 \cong H^0(C_0, \varphi_0^*\mathcal K_X)
\]
 is surjective.
For example, 
 when the canonical sheaf $\mathcal K_{X_0}$ is trivial, then it is clear that this map
 is surjective and also $\varphi_0(C_0)$ is semiregular in the classical sense.
In fact, in this case it is not necessary to assume that the image $\varphi_0(C_0)$ is nodal or reduced,
 and any immersion $\varphi_0$ from a reduced curve $C_0$ is semiregular.  
It is known that when $X_0$ is a K3 surface and the image $\varphi_0(C_0)$
 is reduced, then the map $\varphi_0$
 deforms to general fibers if the class $[\varphi_0(C_0)]$ remains Hodge.
This claim is proved using the twistor family associated with the hyperk\"ahler structure of 
 K3 surfaces, see for example \cite{CGL}.
Corollary \ref{cor:5} gives a generalization of this fact to general surfaces.\\

In \cite{N6}, we also proved the following.
\begin{thm}\label{thm:2}
	Let $X$ be a smooth complex projective surface with an effective canonical class.
	Let $L$ be a very ample class.
	Then, there is a positive number $A$ which depends on $L$ such that 
	for any positive integer $m$, the numerical class of $mL$ contains
	an embedded irreducible nodal curve $C$ whose geometric genus is less than $Am$.\qed
\end{thm}
Such a curve has very large number of nodes, which is roughly $\frac{L^2}{2}m^2$.
In particular, for any positive number $\varepsilon$, we can assume
 $\delta(C) > g(C)^{2-\varepsilon}$ for large $m$, where $\delta(C)$ is the number of nodes of $C$
 and $g(C)$ is the geometric genus of $C$
Moreover, the proof in \cite{N6} shows that 
 we can take $C$ to be semiregular when it is considered as a map $\varphi\colon \tilde C\to X$
 from the normalization $\tilde C$, and $\varphi$
 has non-trivial deformations on which gives
 equisingular deformations of $C$ in $X$.

Now let $Y$ be any smooth projective variety of dimension not less than two.
Let $X$ be a smooth surface in $Y$ which is a complete intersection of sufficiently high degree.
Then by the above theorem there is an embedded irreducible nodal curve $C$
 whose geometric genus is less than $Am$ and which deforms equisingularly in $X$.
Moreover, by Theorem \ref{thm:1}, if we deform $X$ to a smooth surface $X'$ inside $Y$, then the curve $C$
 also deforms to $X'$ equisingularly. 
Since a dense open subset of $Y$ is swept by such surfaces as $X'$, 
 we have the following, which claims that any projective variety can be dominated by 
 an equisingular family of nodal curves which has large number of nodes (and small number of genera). 
\begin{cor}\label{cor:20}
Let $Y$ be a projective variety of dimension $n\geq 2$.
Then for any positive number $\varepsilon$, 
 there is an $(n-1)$-dimensional family $\mathcal C\to B$ of irreducible nodal curves
 whose fibers satisfy $\delta > g^{2-\varepsilon}$, and a map 
 $p\colon \mathcal C\to Y$ which dominates $Y$.
Here $\delta$ is the number of nodes of a fiber of $\mathcal C$ and 
 $g$ is the geometric genus of it.
\end{cor}
\proof
This follows from the above argument applied to a desingularization of $Y$.\qed\\

\section*{Acknowledgement}
\noindent
The author was supported by JSPS KAKENHI Grant Number 18K03313.

\end{document}